\documentclass[12pt]{article}
\usepackage[cp866]{inputenc}
\usepackage[english]{babel}
\usepackage{graphicx}
\usepackage{subfig}
\usepackage{extsizes}
\usepackage{epsfig}
\usepackage[numbers,sort&compress]{natbib}
\usepackage{citehack}
\usepackage{float}
\usepackage{amsmath}
\usepackage{amssymb}
\usepackage{multicol}
\usepackage{amscd}
\usepackage{verbatim}
\usepackage{booktabs}
\usepackage{enumitem}
\newcommand{\vte}{\vartheta}

\newcommand{\p}{\partial}

\makeatletter \@addtoreset{equation}{section} \makeatother
\allowdisplaybreaks[1]

\begin{document}
      \title{On a  kinematic proof of Andoyer variables canonicity}
      \date{}
      \author{Anatoly Neishtadt} 
      \maketitle
    \begin{abstract}
    We present a kinematic proof that  the Andoyer variables in rigid body dynamics are canonical. This  proof is based on the approach of  ``virtual rotations'' by H. Andoyer.  The difference from  the original proof by  Andoyer is  that we do not assume that the fixed in body frame is the frame of principal moments of inertia, and do not use explicit formulas for the kinetic energy of the body that include  moments of inertia. 
       \end{abstract}

\section {Introduction}       
       The canonical Andoyer  variables in rigid body dynamics  were introduced by H. Andoyer in 
       \cite{andoyer_2015,andoyer_2023}.  Wide use of these variables  started after works of V.V. Beletskii \cite{beletskii}, who invented closely related variables, and A. Deprit \cite{deprit}, who invented the same variables as  Andoyer  and used them to represent the free rotation of rigid body in the phase plane. Andoyer proved that his variables are canonical using a kinematic approach  based on ``virtual rotations''.   Deprit used formulas of spherical trigonometry in his proof that the variables are canonical.  Proof of canonicity of Andoyer variables in \cite{deprit_elipe} is based on the representation  of  infinitesimal rotation from the absolute frame of reference to the fixed in body frame via the Andoyer angles. In this note we present a proof based on the original approach by Andoyer. Unlike \cite{andoyer_2015,andoyer_2023, deprit} we do not assume  that the fixed in body frame is  the frame of principal moments of inertia and do not use  explicit formulas for the kinetic energy of the body that include  moments of inertia. We do not calculate the infinitesimal rotation as in \cite{deprit_elipe} either. Our proof uses just general definitions of canonical transformations, kinetic energy and angular momentum as well as properties of the mixed product of vectors. 
       
       
      A review of works on canonical variables in rigid body dynamics prior to Andoyer's studies (F. J. Richelot (1850), J. A. Serret (1866), R. Radau (1869), F. Tisserand (1889)) is given in \cite{deprit_elipe, efr}.

 \section{Coordinate frames. Definition of the Andoyer variables }
\begin{figure}
 \begin{center}
            \includegraphics[scale=0.5, angle=0.0]{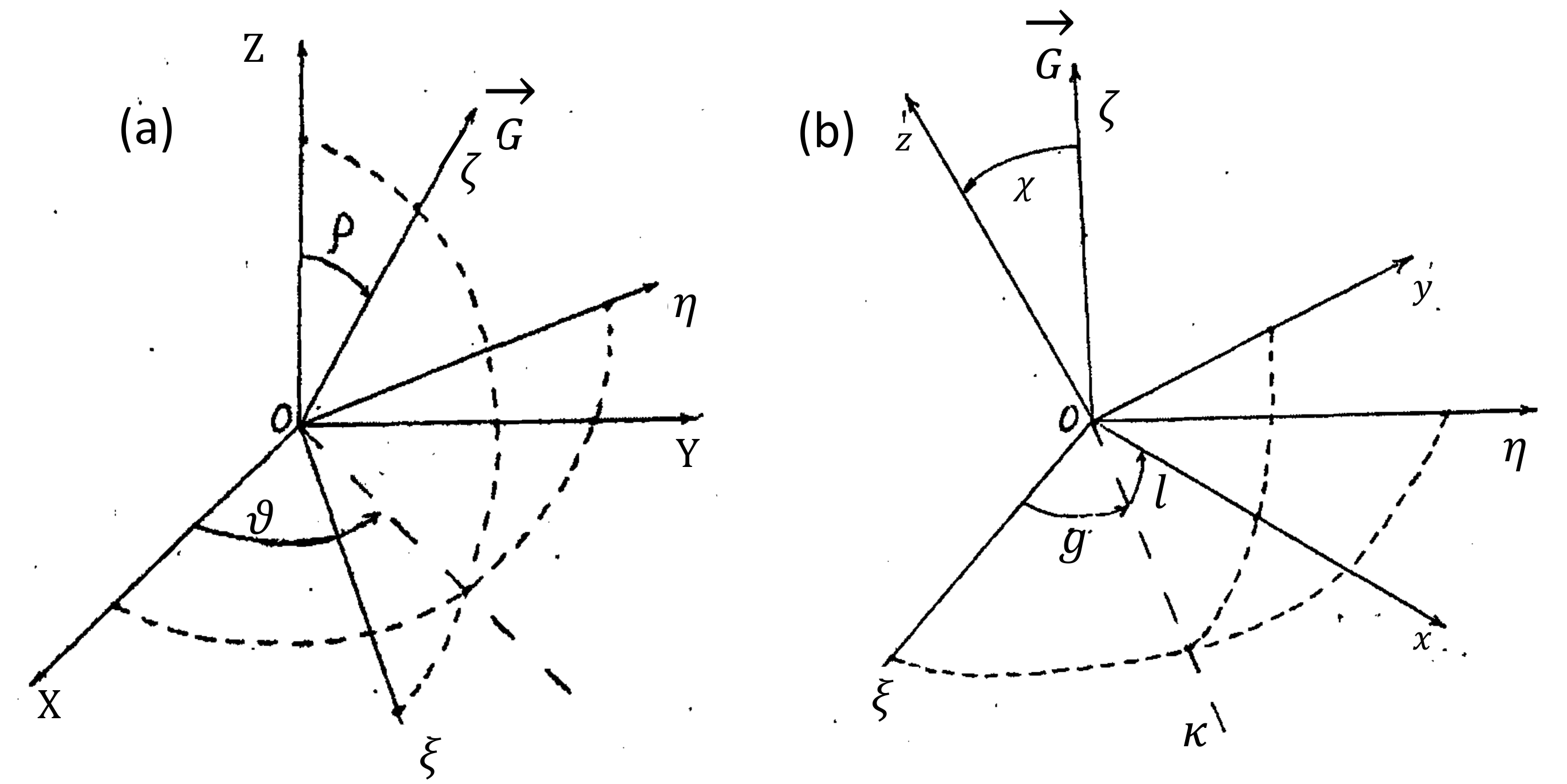}
            \end{center}
           
                      \caption{Coordinate frames.}
            \label{coordinates} 
\end{figure}

 Consider the classical problem of motion of  a rigid body about a fixed point (e.g., \cite{arn_1, goldstein}). Let $OXYZ$ be an absolute Cartesian frame of references, Fig.~ \ref{coordinates},~a. Let $Oxyz$ be a fixed in body Cartesian frame (e.g., the frame of principal moments of inertia of the body for the point $O$  as in \cite{andoyer_2015,andoyer_2023,deprit}), Fig. \ref{coordinates}, b. Denote $\vec G$ the angular momentum of the body. Define a Cartesian frame of reference $O\xi\eta\zeta$ as follows, Fig.1, a, b.   Axis $O\zeta$  is directed along $\vec G$.  Axis $O\xi$ is in the plane $OZ\zeta$. Axis $O\eta$ is orthogonal to the plane $O\xi\zeta$. The canonical Andoyer variables are  $L,G,\Theta, l,g,\vte$ (we a use slightly modified notation from \cite{deprit}). Angles   $l,g,\vte$ are shown in Fig. \ref{coordinates}, a, b, $L$ is the projection of  $\vec G$ onto axis $Oz$, $G$ is the absolute value of $\vec G$, $\Theta$ is the projection of  $\vec G$ onto axis $OZ$.
 Denote   $\vec e_X,\vec e_Y, \vec e_Z$,  $\vec e_x,\vec e_y, \vec e_z$ and  $\vec e_1,\vec e_2, \vec e_3$ unit coordinate vectors of frames $OXYZ$, $Oxyz$ and   $O\xi\eta\zeta$, respectively. We denote  $\vec a\cdot \vec b$ and $\vec a\times \vec b$
 the scalar (``dot'') and the vector (``cross'')  products  of $\vec a$ and $\vec b$.

 \section{Canonicity condition}
 
 Let $q_1,q_2, q_3$ be any generalised coordinates that characterise position of the frame  $Oxyz$ with respect to the absolute frame  $OXYZ$ (e.g., $q_1,q_2, q_3$ could be the Euler angles). The kinetic energy of the body $T$ is a function of these coordinates and their velocities: 
 $T=T( q_1,q_2, q_3, \dot q_1,\dot q_2, \dot q_3)$. Denote  $p_1,p_2, p_3$ the momenta canonically conjugate to these coordinates, $p_j=\p T/\p \dot q_j, \  j=1,2,3$. Then  $q_1,q_2, q_3,p_1,p_2, p_3$ is a system of canonical variables for the considered problem. Express the Andoyer variables via $q_1,q_2, q_3,p_1,p_2, p_3$. To prove that this transformation is  canonical we  should check that (e.g., \cite{arn_1}, p. 241)
 \begin{equation}
 \label{canon_1}
 p_1dq_1+p_2dq_2+p_3dq_3= Ldl+Gdg+\Theta d \vte.
 \end{equation}
 This would prove that the Andoyer variables are canonical.
 
 \medskip

 As it is traditional in Analytical Dynamics, consider the rigid body as a system of $N$ material points. Let $m_i, \vec r_i,\ i=1,2,\ldots, N $ be masses and position vectors  in the absolute frame of these points. Express  position vectors  via coordinates $q_1,q_2, q_3$: $\overline r_i=\overline r_i(q_1,q_2, q_3)$. Then
 \begin{equation}
 \begin{aligned}
  &d {\overline r_i}=\sum_{j=1}^3\frac{\p {\overline r_i}}{\p  q_j}d q_j,\\
 &\dot{\overline r_i}=\sum_{j=1}^3\frac{\p {\overline r_i}}{\p  q_j}\dot q_j .
 \end{aligned}
         \end{equation}
 and
 $$
 \frac{\p \dot{\overline r_i}}{\p \dot q_j}= \frac{\p {\overline r_i}}{\p  q_j}
 $$
 (this is a standard relation used in derivation of the Lagrange equations from the D'Alembert principle, e.g., \cite{goldstein}, p. 20). 
  
    The angular momentum $\vec G$ and the kinetic energy of the body $T$  are
 $$
\vec G= \sum_{i=1}^N m_i (\vec r_i \times \dot {\overline r_i}), \ T=\frac{1}{2}\sum_{i=1}^N m_i (\dot {\overline r_i}\cdot \dot {\overline r_i}).
$$ 
 Consider identities 
 \begin{equation}
 \begin{aligned}
&\sum_{i=1}^N m_i \dot {\overline r_i}\cdot d  {\overline r_i}=\sum_{i=1}^N m_i \dot {\overline r_i}\cdot \sum_{j=1}^3\frac{\p {\overline r_i}}{\p  q_j}d q_j=
\sum_{i=1}^N m_{i} \dot {\overline r_i}\cdot \sum_{j=1}^3\frac{\p \dot{\overline r_i}}{\p \dot q_j}d q_j\\
&=\sum_{j=1}^3\left(\sum_{i=1}^N m_i \dot {\overline r_i}\cdot \frac{\p \dot{\overline r_i}}{\p \dot q_j}
\right)d q_j 
=\sum_{j=1}^3\frac{\p}{\p \dot q_j } \left(\frac{1}{2}\sum_{i=1}^N m_i
 (\dot {\overline r_i}\cdot \dot {\overline r_i })
\right)d q_j \\
&=\sum_{j=1}^3 \frac{\p T}{\p \dot q_j}d q_j= p_1dq_1+p_2dq_2+p_3dq_3.
  \end{aligned}
         \end{equation}
 In view of (\ref{canon_1}) this implies that to prove the canonicity of the Andoyer variables we have to show that
 \begin{equation}
 \label{canon_2}
\sum_{i=1}^N m_i \dot {\overline r_i}\cdot d  {\overline r_i} =Ldl+Gdg+\Theta d \vte.
\end{equation}   
 
\section{The check of canonicity condition (\ref{canon_2})   }

  Express position vectors of  material points via the Andoyer variables. Note that $\overline r_i$ depends on $L, G,\Theta$ via angles $\chi, \rho$, $L=G\cos\chi, \Theta=G\cos\rho$. Thus  $\overline r_i=\overline r_i (l,g,\vte, \chi, \rho), \ i=1,2. \ldots, N$.  Calculate $ d  {\overline r_i} $ and substitute the obtained  expressions into the left hand side of (\ref {canon_2}). We get
 \begin{equation}
 \label{canon_3}
\sum_{i=1}^N m_i \dot {\overline r_i}\cdot d  {\overline r_i} =k_ldl+k_gdg+k_{\vte} d \vte+  k_{\chi}d{\chi}+k_{\rho}d\rho
\end{equation}
with some coefficients $k_l, k_g, \ldots, k_{\rho}$. Calculate these coefficients.

\begin{itemize}
\item
Put    $dl=1$ and  give value 0 to  all other differentials  in the right hand side of    (\ref {canon_3}). Then $d\vec r_i$ should be equal to the velocity of the point with the position vector $\vec r_i$ when the rigid body rotates (counterclockwise) about the axis 
$Oz$ with the angular speed 1: $d \vec r_i=\vec e_z\times \vec r_i$. Then
 \begin{equation*}
 k_l=\sum_{i=1}^N m_i \dot {\overline r_i}\cdot (\vec e_z\times \vec r_i)= (\sum_{i=1}^N m_i \vec r_i \times \dot {\overline r_i})\cdot \vec e_z=\vec G\cdot \vec e_z=L
 \end{equation*}
because  $L$ is the projection of $\vec G$ onto the axis $Oz$.

\item 
Put  $dg=1$ and give value 0 to  all other differentials in the right hand side of    (\ref {canon_3}).  Then $d\vec r_i$ should be equal to the velocity of the point with the position vector $\vec r_i$ when the rigid body rotates   about the axis 
$O\zeta$ with the angular speed 1: $d \vec r_i=\vec e_3\times \vec r_i$. Then
 \begin{equation*}
 k_g=\sum_{i=1}^N m_i \dot {\overline r_i}\cdot (\vec e_3\times \vec r_i)= (\sum_{i=1}^N m_i \vec r_i \times \dot {\overline r_i})\cdot \vec e_3=\vec G\cdot \vec e_3=G
 \end{equation*}
because  $\vec e_3$ is directed along $\vec G$.

\item
Put  $d\vte=1$ and give value 0 to  all other differentials in the right hand side of    (\ref {canon_3}). Then $d\vec r_i$ should be equal to the velocity of the point with the position vector $\vec r_i$ when the rigid body rotates   about the axis 
$OZ$ with the angular speed 1: $d \vec r_i=\vec e_Z\times \vec r_i$. Then
 \begin{equation*}
 k_{\vte}=\sum_{i=1}^N m_i \dot {\overline r_i}\cdot (\vec e_Z\times \vec r_i)= (\sum_{i=1}^N m_i \vec r_i \times \dot {\overline r_i})\cdot \vec e_Z=\vec G\cdot \vec e_Z=\Theta
 \end{equation*}
because  $\Theta$ is the projection of $\vec G$ onto the axis $OZ$.

\item
Put  $d \chi=1$ and give value 0 to  all other differentials in the right hand side of   (\ref {canon_3}). Then $d\vec r_i$ should be equal to the velocity of the point with the position vector $\vec r_i$ when the rigid body rotates   about the node line
$O\kappa$ in Fig.  \ref{coordinates}, b  with the angular speed 1: $d \vec r_i=\vec n\times \vec r_i$, where $\vec n$ is the unit vector of the axis $O\kappa$. Then
 \begin{equation*}
 k_{\chi}=\sum_{i=1}^N m_i \dot {\overline r_i}\cdot (\vec n \times \vec r_i)= (\sum_{i=1}^N m_i \vec r_i \times \dot {\overline r_i})\cdot \vec n=\vec G\cdot \vec n=0
 \end{equation*}
because  $\vec n$ is orthogonal to $\vec G$.

\item
Put  $d \rho=1$ and give value 0 to  all other differentials in the right hand side of    (\ref {canon_3}). Then $d\vec r_i$ should be equal to the velocity of the point with the position vector $\vec r_i$ when the rigid body rotates   about the axis
$O\eta$ with the angular speed 1: $d \vec r_i=\vec 
e_2\times \vec r_i$. Then
 \begin{equation*}
 k_{\rho}=\sum_{i=1}^N m_i \dot {\overline r_i}\cdot (\vec e_2 \times \vec r_i)= (\sum_{i=1}^N m_i \vec r_i \times \dot {\overline r_i})\cdot \vec e_2=\vec G\cdot \vec e_2=0
 \end{equation*}
because  $\vec e_2$ is orthogonal to $\vec G$.
 
 \end{itemize}
 
 Thus we have $k_l=L, k_g=G, k_{\vte}=\Theta, k_{\chi}=k_{\rho}=0$. This proves the canonicity condition  (\ref{canon_2}).
 
 \section{Conlusion}
A kinematic proof is given that  the Andoyer  variables in rigid body dynamics are canonical. The proof is based on  the
approach  of ``virtual rotations'' by Andoyer.  The difference from the original proof by Andoyer is that we do not assume that the fixed in body frame is the frame of principal moments of inertia, and do not use explicit formulas for the kinetic energy of the body that include moments of inertia. 

\vskip 5mm

\noindent Anatoly Neishtadt

\noindent {\small Department of Mathematical Sciences}

\noindent {\small Loughborough University, Loughborough LE11 3TU, United Kingdom}

\noindent {\small Space Research Institute, Moscow 117997, Russia}

\noindent {\footnotesize{E-mail: a.neishtadt@lboro.ac.uk}}

\end{document}